\documentclass[runningheads]{svjour2}
\smartqed  
\usepackage{graphicx}
\usepackage{latexsym}
\usepackage{amssymb,amsmath}
 \DeclareMathOperator{\rk}{rk}
 \DeclareMathOperator{\Ker}{Ker}
 \DeclareMathOperator{\im}{Im}
 \DeclareMathOperator{\Coker}{Coker}
 
 \DeclareMathOperator{\Osc}{Osc}
 \DeclareMathOperator{\Pic}{Pic}

\def\P{\mathcal P}
\def\Q{\mathcal Q}
\def\M{\mathcal M}
\def\L{\mathcal L}
\def\E{\mathcal E}
\def\F{\mathcal F}
\def\N{\mathcal N}

\journalname{Mathematische Zeitschrift}
\begin{document}

\title{Inflectional loci of scrolls
}


\author{Antonio Lanteri \and
        Raquel Mallavibarrena \and
        Ragni Piene}

\authorrunning{A. Lanteri \and R. Mallavibarrena \and R. Piene} 
\institute{A. Lanteri \at
              Dipartimento di Matematica ``F. Enriques'', Universit\`a  Via C.
Saldini, 50, I-20133 Milano, Italy\\
              \email{lanteri@mat.unimi.it}           
           \and
           R. Mallavibarrena \at
              Departamento de Algebra, Facultad de Matem\'aticas, Universidad
Complutense de Madrid, Ciudad Universitaria, E-28040 Madrid, Spain\\
\email{rmallavi@mat.ucm.es}
\and R. Piene \at CMA and Department of Mathematics, University of Oslo, PO Box
1053, Blindern, NO-0316 Oslo, Norway\\
\email{ragnip@math.uio.no} }

\date{Received: date / Revised: date}

\maketitle

\begin{abstract} Let $X\subset \mathbb P^N$ be a scroll over a smooth curve $C$
and let $\L=\mathcal O_{\mathbb P^N}(1)|_X$ denote the hyperplane bundle. The
special geometry of $X$ implies that certain sheaves related to the principal
part bundles of $\L$ are locally free. The inflectional loci of $X$
can be expressed in terms of these sheaves, leading to explicit
formulas for the cohomology classes of the loci. The
formulas imply that the only uninflected scrolls are the balanced rational normal
scrolls.

\end{abstract}

\section{Introduction}
\label{intro} 
Let $X \subset \mathbb P^N$ be a smooth, nondegenerate complex
projective variety of dimension $n$, let $\L=\mathcal O_{\mathbb P^N}(1)|_X$ be
the hyperplane bundle on $X$, and let $V$ be the vector subspace of $H^0(X,\L)$
giving rise to the embedding. Set $V_X = V \otimes \mathcal O_X$. For every
integer $k \geq 0$, let $\P^k_X(\L)$ denote the $k$-th principal part bundle (or
$k$-th jet bundle) of $\L$, and let $j_k:V_X \to \P^k_X(\L)$ be the sheaf
homomorphism sending a section $s \in V$ to its $k$-th jet $j_{k,x}(s)$
evaluated at $x$, for every $x \in X$. Recall that $j_{k,x}(s)$ is represented
in local coordinates by the Taylor expansion of $s$ at $x$, truncated after the
order $k$. The homomorphisms $j_{k,x}$ allow us to define the osculating spaces
to $X$ at $x$ as follows. The $k$-th osculating space to $X$ at $x$ is
$\Osc^k_x(X):=\mathbb P(\im(j_{k,x}))$. Identifying $\mathbb P^N$ with $\mathbb
P(V)$ (the set of codimension $1$ vector subspaces of $V$) we see that
$\Osc^k_x(X)$ is a linear subspace of $\mathbb P^N$.  Since the rank of the
$k$-th jet bundle is $\binom{k+n}n$, we have $\dim\Osc^k_x(X)\leq \binom{k+n}n
-1$.
\par
In the case that $X$ is a scroll over a smooth curve $C$, for $k\geq 2$ this
inequality is strict, for every point $x \in X$. Indeed, let $\pi:X \to C$
denote the structure map. Then, around every point $x \in X$, there are local
coordinates $u,v_2, \ldots , v_n$ such that $u$ is mapped isomorphically by
$\pi$ to a local coordinate on $C$, while $v_2, \ldots , v_n$ are local
coordinates on the fibre through $x$, and every $s \in V$ can be written locally
as $s = a(u)+ \sum_{j=2}^n v_jb_j(u)$, where $a, b_2,
\ldots, b_n$ are regular functions of $u$. Hence, starting with $h=2$, all
derivatives of order $h$ of $s$ vanish, except perhaps $s_{u,u}$, $s_{u,v_j}$,
$j=2, \ldots, n$, for $h=2$, and $s_{u, \ldots ,u,u}$, $s_{u, \ldots, u,v_j}$,
$j=2, \ldots, n$, for $3 \leq h \leq k$. This implies that the rank of
$j_{k,x}$ cannot exceed $kn+1$, hence $\dim\Osc^k_x(X)\leq kn$, at every point
$x \in X$. If equality holds at every point, we say that $X$ is \emph{uninflected}.

Let $X \subset \Bbb P^N$ be a scroll over $C$, and assume $kn \leq N$ and that
the generic rank of $j_k$ is $kn+1$. A main result of this paper is that
the dual $\Q_{k}^{\vee}$ of the sheaf $\Q_k:= \text{Coker}\,j_k$, and the
quotient
sheaf $\E_k:= \P_{X}^k(\L)^{\vee}/\Q_k^{\vee}$, are locally free (Theorem
\ref{thm}). The \emph{$k$-th inflectional locus} $\Phi_{k}$ of $X$ is the set of
points $x\in X$ such that $\rk j_{k,x}< kn+1$.
Now let $k$ be the largest integer such that $kn \leq N$ and
assume
that $\Phi_k$ has the expected codimension $N+1-kn$. Then
$\Phi_k$ has a natural structure as a Cohen--Macaulay scheme, and its cohomology
class is a Segre class of $\E_k$ (Theorem \ref{thm2}). The Segre class can be
computed explicitly (Theorem \ref{thm3}) in terms of Chern classes on the curve
$C$ and the hyperplane bundle class on $X$. In particular, when there are only
finitely many inflection points, their weighted number is equal to
\[
\deg \Phi_{k}=(k+1)\bigl(d+nk(g-1)\bigr).
\]
This formula, and the corresponding formulas in the cases that the expected
dimension of $\Phi_{k}$ is positive, 
allow us to conclude that 
the only  uninflected scrolls are the balanced rational normal scrolls.
\medskip

 In a previous version of this work, we  conjectured 
the degree formula in Corollary \ref{cor1}, but could only show it in special
cases. We
would like to express our deep gratitude to the referee, 
who showed us how Theorem \ref{thm3}, and thus Corollaries \ref{cor1} and
\ref{cor2},
follow from our exact sequences.

\section{The theorems} We want to investigate the inflectional loci of scrolls,
i.e., the loci where the rank of $j_{k,x}$ is smaller than expected. Suppose
that at some point $x \in X$ (hence at a general point), we have
$\rk(j_{k,x})=kn+1$, and let $\Q_k$ denote the cokernel of the map $j_k$. Then
we have an exact sequence of sheaves
\begin{equation}
\label{1}
 0 \to \mathcal P_k \to \P^k_X(\L) \to \Q_k \to 0,
\end{equation}
where $\mathcal P_k := \im(j_k)$.

The $k$-th inflectional locus $\Phi_k$ of $X$ is the locus where the map
$j_k:V_X \to \P^k_X(\L)$ does not have the maximal rank $kn+1$, hence where the
sheaf $\mathcal P_k$ is not a vector bundle. It is also the locus where the dual
map $j_k^{\vee}:\P^k_X(\L)^{\vee}\to V^{\vee}_X$ does not have maximal rank.
Clearly, for $k \geq k'$ we have
 $\Phi_k \supseteq \Phi_{k'}$, because of the surjections
 $\P_k\to \P_k'$.
 
 For simplicity, we shall denote by $\Omega_{Y}$ the cotangent sheaf
$\Omega_{Y}^1$ of a variety $Y$. We shall write $T_{Y}$ for the tangent sheaf,
equal to the dual sheaf $\Omega_{Y}^\vee$. If $\N$ is a line bundle, we write
$\N^{i}$ instead of $\N^{\otimes i}$ and $\N^{-1}$ for $\N^\vee$.

\begin{theorem}\label{thm} 
Let $X \subset \mathbb P^N$ be a $n$-dimensional
scroll over a smooth curve $C$, with hyperplane bundle $\L=\mathcal O_{\mathbb
P^N}(1)|_{X}$. For all $k \geq 1$ such that $kn\le N$, assume the generic rank
of $j_{k}$ is $kn+1$, and set $\Q_k=\Coker j_{k}$ as above. Then, for such $k$,
\begin{itemize}
\item[\textrm{(i)}] $\Q_k^{\vee}$ is a locally free sheaf of rank $\binom{n+k}n
- (kn+1)$, \item[\textrm{(ii)}] there exist locally free sheaves $\M_k$ of rank
$\binom{n+k-1}{n-1}-n$ and exact sequences
\begin{equation}
\label{th1}
0 \to \Q_{k-1}^{\vee} \to  \Q_{k}^{\vee} \to \M_k^{\vee}
\to 0,
\end{equation}
and
\begin{equation}
\label{th2}
 0 \to \pi^*\Omega^{k-1}_C \otimes
\Omega_X \otimes \L \to S^k \Omega_X \otimes \L \to \M_k \to 0,
\end{equation}
\item[\textrm{(iii)}] the quotient sheaves
$\E_{k}:=\P^k_X(\L)^{\vee}/\Q_k^{\vee}$ are locally free, of rank $kn+1$, and
there exist exact sequences
\begin{equation}
\label{th3} 
0 \to \E_{k-1} \to \E_{k} \to \pi^*T^{k-1}_C
\otimes T_X \otimes \L^{-1} \to 0.
\end{equation}
\end{itemize}
\end{theorem}

\begin{proof} There is an obvious inclusion
\[
\pi^*\Omega^{k-1}_C \otimes \Omega_X \otimes \L
\subset S^{k-1}\Omega_X \otimes \Omega_X \otimes \L.
\]
 By restricting the natural homomorphism $S^{k-1}\Omega_X \otimes
\Omega_X \otimes \L \to S^{k}\Omega_X \otimes \L$ we get an injective and
locally split homomorphism
\begin{equation}
\label{iota}
\iota_k: \pi^*\Omega^{k-1}_C \otimes \Omega_X \otimes \L
\to S^{k}\Omega_X \otimes \L.
\end{equation}
To see this, let $x \in X$ and let $u$ denote a local coordinate on the base
curve $C$ around $\pi(x)$ and $v_2, \ldots , v_n$ local coordinates in the fibre
of $X$ through $x$, around $x$. Then $u, v_2, \ldots ,v_n$ are local coordinates
on $X$ around $x$. So, letting $A:= \mathcal O_{X,x}$, we have the following
isomorphisms: $$
\textstyle (\pi^*\Omega_C)_x \cong A\ du, \qquad \Omega_{X,x} \cong A\ du \oplus
\bigoplus_{i=2}^{n} A\ dv_i,$$
and $$
\textstyle (S^{k-1}\Omega_X)_x \cong \bigoplus_i A\ du^{i_1} dv_2^{i_2} \dots
dv_n^{i_n}, \quad \textrm{ with } \quad \sum_j i_j = k-1.$$
The map $\iota_{k,x}$ on the stalks is clearly injective, since it acts on the
differentials in the following way: $$
du^{k-1} \otimes du \mapsto du^k \qquad \textrm{and} \qquad du^{k-1} \otimes
dv_i \mapsto du^{k-1} dv_i.$$
The same is true for the map $\iota_{k}(x)$  on the fibres. Hence $\iota_{k}$ is
locally split. Set $\M_k:= \Coker
\iota_k$. It follows that $\M_k$ is locally free, and we get the second exact
sequence (\ref{th2}).

Let $\E_{k}:=\P^k_X(\L)^{\vee}/\Q_k^{\vee}$ denote the quotient sheaf, and
 consider the following diagram, all of whose horizontal sequences are exact: $$
\begin{matrix} &  & 0 &  & 0 &  & 0 & & \\
 &  & \downarrow &  & \downarrow &  & \downarrow &  &\\
 0 & \to & \Q_{k-1}^{\vee} & \to & \mathcal P^{k-1}_{X}(\L)^{\vee} & \to &
\E_{k-1} & \to & 0\\
 &  & \downarrow &  & \downarrow &  & \downarrow & & \\
 0 & \to & \Q_k^{\vee} & \to & \P^k_X(\L)^{\vee} & \to & \E_k & \to & 0\\
 &  &  &  & \downarrow &  &  &  & \\
 0 & \to & \M_k^{\vee} & \to & (S^k\Omega_X \otimes \L)^{\vee} & \overset
{\iota_{k}^{\vee}} \longrightarrow&
 (\pi^*\Omega^{k-1}_C \otimes \Omega_X \otimes \L)^{\vee} & \to & 0\\
 &  &  &  & \downarrow &  &  & &  \\
 &  &  &  &  0  & & & &\\
\end{matrix}$$
We will complete it to a commutative diagram in which all vertical sequences are
also exact.

Let us first consider the composition $\alpha$ of the map $\iota_{k}$
(\ref{iota}) with the maps $S^k\Omega_X \otimes \L \to \P^k_X(\L)$ and
$\P^k_X(\L)\to \Q_k$,
\[
\alpha\colon \pi^*\Omega^{k-1}_C \otimes \Omega_X \otimes \L \to  \Q_k.
\]
We want to show that this map is generically $0$. Let $x \in X$ be a general
point and use local coordinates as above. Then $(\P^k_X(\L))_x \cong \mathcal
P^{k-1}_{X}(\L)_x \oplus (S^k\Omega_X\otimes \L)_x$ and $\alpha_x$ sends the
generators $du^{k-1} \otimes du$ and $du^{k-1}\otimes dv_i$ to $du^{k}$ and
$du^{k-1}dv_i$ in the second summand. But these elements are in
$\im(j_{k,x})=(\P_{k})_{x}$, hence they go to $0$ in $(\Q_{k})_{x}$, by
(\ref{1}).

Since $\alpha$ is generically $0$, so is its dual,
\[
\alpha^{\vee}\colon \Q_k^{\vee} \to (\pi^*\Omega^{k-1}_C \otimes \Omega_X
\otimes
\L)^{\vee}=\pi^*T_{C}^{k-1}\otimes T_{X}\otimes \L^{-1}.
\]
Since the target sheaf is locally free, it has no torsion, hence $\alpha^{\vee}$
is everywhere zero. Therefore we get induced maps $\psi$ and $\beta$ making the
diagram commute: $$
\begin{matrix}  &  & 0 &  & 0 &  & 0 & & \\
 &  & \downarrow &  & \downarrow &  & \downarrow &  &\\
 0 & \to & \Q_{k-1}^{\vee} & \to & \mathcal P^{k-1}_{X}(\L)^{\vee} & \to &
\E_{k-1} & \to & 0\\
 &  & \downarrow &  & \downarrow &  & \downarrow & & \\
 0 & \to & \Q_k^{\vee} & \to & \P^k_X(\L)^{\vee} & \to & \E_k & \to & 0\\
 &  & \psi \downarrow &  & \downarrow &  & \beta \downarrow &  & \\
 0 & \to & \M_k^{\vee} & \to & (S^k\Omega_X \otimes \L)^{\vee} & \to &
 \pi^*T_{C}^{k-1}\otimes T_{X}\otimes \L^{-1} & \to & 0\\
 &  & \downarrow  &  & \downarrow &  & \downarrow & &  \\
 &  & 0 &  &  0  & & 0 & &.\\
\end{matrix}$$

We show that: a) $\psi$ is surjective, and b) $\Ker \psi =
\mathcal Q_{k-1}^{\vee}$. First note that $\E_{ k-1}\to \E_{k}$ is injective,
since both sheaves are subsheaves of $V_{X}^{\vee}$. Then fact a) follows from
the snake lemma: we have the exact sequence $\mathcal P^{k-1}_{X}(\L)^{\vee} \to
\mathcal \E_{k-1}
\to \Coker \psi \to 0 \to 0$; the first map is surjective by the diagram, hence
the second is zero; thus the third is injective, but since its image is zero, we
conclude that $\Coker \psi$ itself is zero. As for b), clearly $Q_{k-1}^{\vee}
\subseteq \Ker \psi$, by an easy diagram chase. The converse also follows by a
diagram chase: Let $\xi \in (\Ker \psi)_x$ for some $x \in X$. Since $\xi$ goes
to zero in $(\M_k^{\vee})_x$, then its image, say $y$ in
$(\P^k_X(\L))^{\vee}_x$, goes to zero via both the horizontal and the vertical
maps. Then $y$ comes from an element $z \in (\mathcal
P^{k-1}_{X}(\L)^{\vee})_x$, which must go to zero by the horizontal map, due to
the commutativity of the right-upper square. Thus there exists a $w \in
(Q_{k-1}^{\vee})_x$ mapping to $z$. Since the map $\xi \mapsto y$ is injective,
we thus conclude that $\xi$ is the image of $w$: this shows that $(\Ker \psi)_x
\subseteq (Q_{k-1}^{\vee})_x$.

To conclude, for every $k \geq 2$ the first vertical sequence gives the exact
sequence (\ref{th1}). In particular, since $\M_k^{\vee}$ is locally free and
$\Q_2^{\vee}=\M_{2}^{\vee}$ because $\Q_1 = 0$, this shows, by induction, that
$\mathcal Q_k^{\vee}$ is locally free for every $k \geq 2$. The assertion on the
rank of $\Q_k^{\vee}$ follows from the fact that it equals the generic rank of
$\Q_k$, which is given by $\rk\P^k_X(\L) - \rk j_{k,x}$ for a general $x \in X$.
This proves (i) and (ii).

To prove (iii), observe that since $\Q_{1}=0$, $\E_{1}\cong 
\P^{1}_{X}(\L)^{\vee}$ is locally free. Using the exactness of the rightmost
vertical sequence in the commutative diagram with $k=2$, we
 conclude that $\E_{2}$ must be locally free, since an extension of two locally
free sheaves is locally free. Hence we deduce, recursively, that all $\E_{k}$
are locally free. \qed
\end{proof}

\begin{theorem}\label{thm2} Let $X \subset \mathbb P^N$ be a $n$-dimensional
scroll over a smooth curve $C$, with hyperplane bundle $\L=\mathcal O_{\mathbb
P^N}(1)|_{X}$. Let $k$ be the largest integer such that $kn\le N$ and assume
that the generic rank of $j_{k}$ is $kn+1$. If the $k$-th inflectional locus
$\Phi_{k}$  of $X$ has codimension $\ell :=N+1-kn$ or is empty, then it has a
natural structure as a Cohen--Macaulay scheme, and its class is equal to the
$\ell$th term of the Segre class of $\E_{k}$,
\begin{equation*}
[\Phi_{k}]=\bigl[c(\E_{k})^{-1}]_{\ell}\,\,,
\end{equation*}
where $\E_{k}=\P^k_X(\L)^{\vee}/\Q_k^{\vee}$.
\end{theorem}

\begin{proof} Note that the assumptions imply that the generic rank of
$j_{k^{\prime}}$ is $k^{\prime}n+1$ for $k^{\prime}\le k$, so that the
assumptions of Theorem
\ref{thm} are satisfied.

It follows from the definition that the $k$-th inflectional locus is equal to the
degeneracy locus of the map of locally free sheaves
\[
0\to \E_{k}=\P^k_X(\L)^{\vee}/Q_k^{\vee}\to V_{X}^{\vee}\,,
\]
hence it has a natural structure as a Cohen--Macaulay scheme when it has the
expected codimension $N+1-kn$. By Porteous' formula \cite[Ex.~14.4.1,
p.~255]{Fu}, the class of the $k$-th inflectional locus is equal to the class
\begin{equation*}
\bigl[c(V^{\vee}_{X})c(\E_{k})^{-1}\bigr]_{\ell}
=\bigl[c(\E_{k})^{-1}\bigr]_{\ell}\,
\end{equation*}
where $\ell = N+1-kn$. 
\qed
\end{proof}

Note that, since $k$ is the largest integer such that $kn\le N$, we also have
$N\le (k+1)n-1$. This implies that
\begin{equation*}
1 \leq \ell \leq n\,.
\end{equation*}

\begin{theorem}\label{thm3} The $j$-th term of $c(\E_{k})^{-1}$, for
$j=1,\ldots,n$, is equal to
\[
L^j+k\bigl(d+(n(k-1)+2j)(g-1)\bigr)L^{j-1}F,
\]
where $L=c_{1}(\L)$ denotes the class of a hyperplane section of $X$, 
$F$ is the class of a fiber of the map $\pi\colon X\to C$, $d$ is the degree of
$X$, and $g$ is the genus of $C$.
\end{theorem}

\begin{proof} We use the exact sequences (\ref{th3}) for $i=1,\ldots,k$ to get
\[
c(\E_{k})=c(\pi^*T_{C}^{k-1}\otimes T_{X}\otimes
\L^{-1})c(\E_{k-1})=\prod_{i=0}^{k-1}c(\pi^*T_{C}^{i}\otimes T_{X}\otimes
\L^{-1})c(\L^{-1}).
\]
The standard exact sequence
\[
0\to \pi^*\Omega_{C}\to \Omega_{X}\to \Omega_{X/C}\to 0
\]
gives, by dualizing and tensoring with $\pi^*T_{C}^{i}$ and $\L^{-1}$, exact
sequences
\[
0\to \pi^*T_{C}^{i}\otimes T_{X/C}\otimes \L^{-1}\to 
\pi^*T_{C}^{i}\otimes T_{X}\otimes \L^{-1}\to 
\pi^*T_{C}^{i+1}\otimes \L^{-1}\to 0.
\]
 The sheaf $\F:=\pi_{*}\L$ is locally free, with rank $n$, and we have the
standard exact sequence
\[
0\to \Omega_{X/C}\otimes \L\to \pi^*\F\to \L\to 0.
\]
Dualizing and tensoring with $\pi^*T_{C}^{i}$, we obtain sequences
\[
0\to \pi^*T_{C}^{i}\otimes \L^{-1}\to \pi^*(T_{C}^{i}\otimes\F^{\vee} )
\to \pi^*T_{C}^{i}\otimes T_{X/C}\otimes \L^{-1}\to 0.
\]
Hence 
\[
c(\pi^*T_{C}^{i}\otimes T_{X}\otimes\L^{-1})=c(\pi^*T_{C}^{i+1}\otimes
\L^{-1})c( \pi^*(T_{C}^{i}\otimes\F^{\vee}))c(\pi^*T_{C}^{i}\otimes
\L^{-1})^{-1},
\]
which gives, because of cancellations in the product, the expression
\[
c(\E_{k})=\prod_{i=0}^{k-1}\pi^*c(\F^\vee\otimes T_{C}^{i})
c(\pi^*T_{C}^{k}\otimes \L^{-1}).
\]
The last Chern class in this product is the Chern class of a line bundle, so
that $c(\pi^*T_{C}^{k}\otimes
\L^{-1})=1+k\pi^*c_{1}(T_{C})-c_{1}(\L)=1-k(2g-2)F-L$, since
$\pi^*c_{1}(T_{C})=-\pi^*c_{1}(\Omega_{C})=-(2g-2)F$. Since
$\F^\vee\otimes
T_{C}^{i}$ is a bundle on the curve $C$, its Chern class is just
$1+c_{1}(\F^\vee\otimes T_{C}^{i})=1-c_{1}(\pi_{*}\L)+nic_{1}(T_{C})$, and its
inverse Chern class is $1+c_{1}(\pi_{*}\L)+nic_{1}(\Omega_{C})$. Because
$\pi^*c_{1}(\pi_{*}\L)=dF$, we get
$\pi^*(1+c_{1}(\pi_{*}\L)+nic_{1}(\Omega_{C}))=1+(d+2in(g-1))F$, and thus
\begin{eqnarray*}
c(\E_{k})^{-1}&=\prod_{i=0}^{k-1}\bigl(1+(d+2in(g-1))F\bigr)(1-2k(g-1)F-L)^{-1}\\
&=(1+aF)(1-bF-L)^{-1},
\end{eqnarray*}
where we set $a:=k(d+n(k-1)(g-1))$ and $b:=2k(g-1)$ and used the fact that
$F^{i}=0$ for $i>1$.  The $j$-th term of this class is equal to
\[
(bF+L)^{j}+aF (bF+L)^{j-1}=L^{j}+jbL^{j-1}F+aL^{j-1}F=L^{j}+(a+jb)L^{j-1}F,
\]
which is what we wanted to prove.\qed
\end{proof}

\begin{corollary}\label{cor1} Under the assumptions of Theorem \ref{thm2}, the
class of the inflectional locus of $X$ is equal to
\[
[\Phi_{k}]=L^{N+1-kn}+k\bigl(d+(n(k-1)+2(N+1-kn))(g-1)\bigr)L^{N-kn}F,
\]
and its degree is equal to
\[
\deg \Phi_{k}=(k+1)d+k\bigl(2(N+1)-(k+1)n\bigr)(g-1).
\]
In particular, if $N=(k+1)n-1$, then $\Phi_{k}$ is $0$-dimensional, and its
degree is equal to
\begin{equation}\label{deg}
\deg \Phi_{k}=(k+1)\bigl(d+nk(g-1)\bigr).
\end{equation}
\end{corollary}

\begin{corollary}\label{cor2} Let $\ell$ be an integer, $1\le \ell \le n$. The
only uninflected scroll $X\subset
\mathbb P^{kn+\ell-1}$ of dimension $n$ is the balanced rational normal scroll
of degree $kn$ in $\mathbb P^{(k+1)n-1}$.
\end{corollary}
\begin{proof} If $X$ is uninflected, then the assumptions of Theorem
\ref{thm2} are satisfied, since $\Phi_{k}=\emptyset$. By Corollary \ref{cor1}
the class 
\[
L^{\ell}+k\bigl(d+(n(k-1)+2\ell)(g-1)\bigr)L^{\ell-1}F
\]
is $0$. If $\ell < n$, we can intersect this class with $L^{n-\ell-1}F$ and
obtain
$L^{n-1}F=0$, using the fact that $F^2=0$. But
$L^{n-1}F=1$ is the degree of the linear space $F$, thus we get a contradiction. 

Hence we may assume $\ell=n$, so that $N=(k+1)n-1$. Setting $\deg \Phi_{k}=0$
in (\ref{deg}) implies
$g=0$ and $d=kn$. Therefore, $X$ is a smooth, nondegenerate
rational scroll of minimal degree, hence it is linearly normal. The explicit
description of the maps $j_{k}$ given in
\cite[p. 1050]{PS} shows that the only uninflected rational normal
$n$-dimensional scrolls in $\mathbb P^{(k+1)n-1}$are
the balanced ones, i.~e., the ones given by $X=\mathbb P (\pi_{*}\L)=\mathbb
P(\mathcal O_{\mathbb P^1}(k)\oplus\ldots \oplus
\mathcal O_{\mathbb P^1}(k))$ on $C=\mathbb P^{1}$.
\qed
\end{proof}

\section{Examples} 

In this section we give geometric descriptions and details
about inflectional loci in some particular, but relevant, cases.

In the situation of Theorem \ref{thm2}, when  $n=1$, we have
$X=C$ and $N=k$, so that $X\subset \mathbb P^k$ is a nondegenerate curve.
Corollary
\ref{cor1} gives the formula
\begin{equation}\label{n=1}
\deg \Phi_{k}=(k+1)\bigl(d+k(g-1)\bigr)
\end{equation}
for the total (weighted) number of inflection points. This classical
formula, valid also when $X$ is singular, goes back to Veronese and has been
reproved many times (see e. g. \cite[Thm. 3.2]{P}).
\medskip

When $n=2$ and $k=2$, we have a surface scroll $X\subset
\mathbb P^5$.  In this case, Corollaries
\ref{cor1} and \ref{cor2} were shown by Shifrin \cite[Prop. 4.3 and Thm. 4.3,
p. 247]{Sh}; in the more general case of a surface scroll $X\subset \mathbb P^{2k+1}$,
Corollary
\ref{cor2} was shown by Piene and Tai \cite[p. 221]{PT}.
\medskip

Note that for $n=2$ and $X \subset \mathbb P^{2k+2}$, we are outside the
range of Corollary \ref{cor2}. In this case there are several examples of 
uninflected
scrolls: for example, the scroll with $g=1$, defined by an indecomposable rank
$2$ vector bundle of degree $2k+3$  \cite[Thm. A]{L}, and scrolls with
$g=0$, both normal (the semibalanced scroll $\mathbb P(\mathcal O_{\mathbb
P^1}(k)
\oplus \mathcal O_{\mathbb P^1}(k+1))$ \cite{PT})  and non-normal \cite[Thm.
3.4]{LM2}.
\medskip

In the next example we consider at the same time the following
cases:
\begin{itemize}
\item[(i)] $g=0$, $k \geq 2$; 
\item[(ii)] $g=1$, $k
\geq 3$.
\end{itemize}
For $i=1, \dots ,n$, consider line bundles $\L_i \in \Pic(C)$ such
that $\deg \L_i=g+k-1$ for $i=1, \dots ,n-1$ and $\deg \L_n =
g+k$. Take $X = \mathbb P(\F)$, where $\F =\oplus_{i=1}^n \L_i$, and
let $\L$ be the tautological line bundle on $X$. Clearly $\L$ is
very ample; moreover, $h^0(\L)=h^0(\F)=(n-1)k+k+1=nk+1$. So, $X$
embedded by $|\L|$ is a linearly normal scroll in $\mathbb P^{kn}$
of degree $d= \deg \F = n(g+k)-(n-1)$ in both cases (i)
and (ii). Let $C_i$ be the generating section of the scroll $X$
corresponding to the $i$-th summand $\L_i$ of $\F$.
Note that $C_i$ is embedded in $\mathbb P^{k-1}$ as a rational
(resp. elliptic) normal curve in case (i) (resp. (ii) ) for
$i=1, \dots , n-1$. The same holds for $C_n$ in $\mathbb P^k$.
Thus $\Phi_{k-1}(C_i) = \emptyset$ and $\Phi_{k}(C_i)=C_i$ for
$i=1, \dots ,n-1$ in case (i), while $\Phi_{k-1}(C_i)$ consists of
$k^2$ points, according to (\ref{n=1}), and $\Phi_{k}(C_i)=C_i$
for $i=1, \dots ,n-1$ in case (ii). Similarly, $\Phi_k(C_n)$ is
either empty (case (i)) or consists of $(k+1)^2$ points (case (ii)).
 Let $Y:=\mathbb P(\oplus_{i=1}^{n-1} \L_i)$ denote the
 $(n-1)$-dimensional sub-scroll of $X$ generated
by the sections $C_i$, for $i=1, \dots , n-1$. The above facts
imply that the $k$-th inflectional locus $\Phi_k$ of $X$ is equal to $Y$ in case
(i), and to
the union of $Y$ and $(n-1)k^2 + (k+1)^2$ fibers in case (ii). This follows
from
\cite[p.\
1050]{PS} in case (i) and \cite[p.\ 152]{MP} in case (ii) (see also
\cite[Prop. 2.6 and Cor. 2.9]{LM2}). We have $[Y]=L-(g+k)F$, since 
$\deg \L_n = g+k$. This gives $[\Phi_{k}]=L-kF$ in case (i), and 
$[\Phi_{k}]=L-(k+1)F+((n-1)k^2 + (k+1)^2)F=L+k(nk+1)F$ in case (ii).
This agrees with
Corollary \ref{cor1}, which gives $[\Phi_k] = L-kF$ when  $N=kn$ and $g=0$, and
$[\Phi_k] =L+kd F$ when $N=kn$ and $g=1$. 
\medskip

Note that the only rational nondegenerate 
scroll  in $\mathbb P^{2n}$ is the linearly normal
rational scroll of degree $n+1$ considered in case (i) above, with $k=2$. In
fact, for any smooth $n$-dimensional scroll  $X\subset\mathbb P^{2n}$,
the well known double point (or self-intersection) formula becomes
 \begin{equation}\label{dbl}
 (d-n)(d-n-1)=n(n+1)g,
 \end{equation}
so that for $g=0$, we must have $d=n+1$ if $X$ is nondegenerate.

It is in fact conjectured that any
scroll $X \subset \mathbb P^{2n}$ of dimension $n$ has
$g=0$ or $1$; this conjecture holds for $n \leq 4$
\cite[Cor. 5]{IT}. 
\medskip

For $g=1$, the case $k=2$ is not covered by (ii) in the above example. Note that,
by
(\ref{dbl}), such a scroll must have degree $2n+1$. In fact, it is well known
that the only such scrolls are the ones constructed as follows.
Consider
a smooth, elliptic curve $C$, and define inductively rank $i$, degree $1$
sheaves $\F_{i}$ by starting with $\F_1 = \mathcal O_C(p)$, for some $p\in C$,
and using the non-split exact sequences $0 \to \mathcal O_C 
\to \F_{i+1} \to \F_{i} \to 0$. Taking  points $p_1, p_2
\in C$, then
$X = \mathbb P(\F_n(p_1+p_2))$ can be embedded by the
tautological line bundle, giving an indecomposable scroll of degree $2n+1$ in
$\mathbb P^{2n}$. If the assumptions of Theorem
\ref{thm2} are
satisfied, then Corollary \ref{cor1} gives $[\Phi_2] = L + 2(2n+1)F$.
\medskip

\begin{acknowledgements} 
The first author would like to thank
the MUR of the Italian Government for support received in the framework
of the PRIN ``Geometry on Algebraic Varieties'' (Cofin 2002 and 2004),
as well as the University of Milan (FIRST) for making this collaboration
possible.           
The second author wants to thank for the funds supporting this research from the
projects BFM2003-03917/MATE (Spanish Ministry of Education) and Santander/UCM
PR27/05-138.        
\end{acknowledgements}


\begin{thebibliography}{ACGH}

\bibitem {Fu} Fulton, W.:
      Intersection Theory. 2nd ed. Springer-Verlag, New York--Heidelberg--Berlin
(1998).

\bibitem {IT} Ionescu, P., Toma, M.: On very ample vector bundles
on curves. Internat. J. Math. {\bf 8}, 633--643 (1997).

\bibitem {L} Lanteri, A.: On the osculatory behavior of surface scrolls. 
Matematiche (Catania) {\bf 55}, 447--458 (2000).

\bibitem {LM2} Lanteri, A., Mallavibarrena, R.:  Osculating
properties of decomposable scrolls, Preprint (2006).

\bibitem {MP} Mallavibarrena, R., Piene, R.: Duality for elliptic normal surface
scrolls. Contemp. Math. {\bf 123}, 149--160 (1991).

\bibitem {P} Piene, R.:
Numerical characters of a curve in
projective $n$-space.  In: Holm, P. (ed.) Real and complex singularities.
Proceedings,
Oslo 1976, pp.~475-496. Sijthoff \& Noordhoff  (1977).

\bibitem {PS} Piene, R., Sacchiero, G.: Duality for rational normal scrolls.
Comm. Algebra {\bf 12}, 1041-1066 (1984).

\bibitem {PT} Piene, R., Tai, H. S.:  A characterization of balanced rational
normal scrolls in terms of their osculating spaces. In Xambo-Descamps, S. (ed.)
Enumerative Geometry, Proc. Sitges, 1987, pp.~ 215--224. Lecture Notes in Math.
{\bf  1436 }, Springer-Verlag  (1990).

\bibitem {Sh} Shifrin, T.: The osculatory behavior of surfaces in $\mathbb P^5$.
 Pacif. J. Math. {\bf 123},
  227--256  (1986).

\end{thebibliography}
\end{document}